\documentclass[12pt,reqno]{amsart}
\usepackage{amssymb,amsmath}
\usepackage[final]{graphicx}
\usepackage{epic}
\usepackage{pstricks}
\usepackage{color}
\usepackage[ansinew]{inputenc}
\usepackage[english]{babel}
\usepackage{latexsym}
\usepackage{eufrak}
\usepackage{mathrsfs}
\usepackage{fancybox}
\usepackage{fancyhdr}
\usepackage{array}
\newtheorem{thm}{Theorem}

\newtheorem{lem}[thm]{Lemma}

\def\C{\mathbb{C}}

\title[linearization coefficients of Bessel polynomials]{An explicit formula for the linearization coefficients of Bessel polynomials II}
\author{Mohamed Jalel Atia, Gabes univ, Tunisia }

\date{\today}
\begin{document}
\maketitle

\begin{abstract}
In this paper, a single sum  formula for  the linearization coefficients  of the Bessel polynomials is given. In three special cases
this formula reduces indeed to
either Atia and Zeng's formula (Ramanujan Journal, Doi 10.1007/s11139-011-9348-4)
or Berg and Vignat's formulas  in their proof of the  positivity results about
these  coefficients  (Constructive Approximation, {\bf 27} (2008), 15-32). As a bonus, a formula reducing a sum of hypergeometric functions $_3F_2$ to $_2F_1$ is obtained.
 \end{abstract}
\medskip


\section{Introduction}

The Bessel polynomials $q_n$ of degree $n$ are defined by
\begin{equation} \label{poly_bess}
q_n(u)=\sum_{k=0}^n{(-n)_k2^k\over (-2n)_kk!}u^k,
\end{equation}
where we use the Pochhammer symbol $(z)_n:=z(z+1)\ldots (z+n-1)$ for $z\in \C$, $n=0,1,...$.
The first values are
$$
q_0(u)=1,\quad q_1(u)=1+u, \quad q_2(u)=1+u+\frac{u^2}{3}.
$$
Some recursion formulas for $q_n$ are
\begin{equation} \label{recu_bess}
q_{n+1}(u)=q_n(u)+{u^2\over 4n^2-1}q_{n-1}(u),\ n\geq 1,
\end{equation}
\begin{equation} \label{recu_bess der}
q_{n}'(u)=q_n(u)-{u\over 2n-1}q_{n-1}(u),\ n\geq 1.
\end{equation}

Using hypergeometric functions, we have $q_n(u)=\ _1F_1(-n;-2n;2u)$.
They are normalized according to $q_n(0)=1$, and thus differ from the monic polynomials $\theta_n(u)$ in Grosswald's monograph \cite{Gro}:
$$\theta_n(u)=\displaystyle{(2n)!\over n!2^n}q_n(u).$$
The polynomials $\theta_n$ are sometimes called the reverse Bessel polynomials and $y_n(u)=u^n\theta_n({1\over u})$ the ordinary Bessel polynomials.
These Bessel polynomials are, then, written as
\begin{equation}\label{yn}
y_n(u)=\displaystyle{(2n)!\over n!2^n}u^nq_n({1\over u})=\sum_{k=0}^n{(n+k)!\over 2^k k!(n-k)!}u^k.
\end{equation}
The {\it linearization problem} is the problem of finding the coefficients $\beta_{k}^{(n,m)}(a_1,a_2)$ in the expansion of the product $P_n(a_1u)Q_m(a_2u)$ of two polynomials systems in terms of a third sequence of polynomials $R_k(u)$,
\begin{equation} \label{eq:genlin}
P_n(a_1u)Q_{m}(a_2u)=\sum_{k=0}^{n+m}\beta_{k}^{(n,m)}(a_1,a_2)R_{k}(u).
\end{equation}

The polynomials $P_n$, $Q_m$ and $R_k$ belong to three different polynomial families.
In the case $P=Q=R$ and $a_1=a_2=1$, we get the (standard) linearization or Clebsch-Gordan-type problem. If $Q_m(u)\equiv 1$, we are faced with the so-called connection problem.\\

In the case $P=Q=R$ and $a_1=a,\ a_2=1-a$, we get the Berg-Vignat linearization problem.
And, finally, in the case $P=Q=R$ and for any $a_1,a_2$, we get a new linearization problem.

In this paper, we are interested by this new linearization problem and by the linearization  coefficients $\beta_{k}^{(n,m)}(a_1,a_2)$  in the case of the Bessel polynomials which are defined by

\begin{equation} \label{eq:lin}
q_n(a_1u)q_{m}(a_2u)=\sum_{k=0}^{n+m}\beta_{k}^{(n,m)}(a_1,a_2)q_{k}(u).
\end{equation}
For example,  we have
\begin{equation}
q_3(a_1u)q_{5}(a_2u)=\sum_{k=0}^{8}\beta_{k}^{(3,5)}(a_1,a_2)q_{k}(u)
\end{equation}
where
\begin{align*}
\beta_8^{(3,5)}(a_1,a_2)&=143\,{a_1}^{3} a_2 ^{5},\\
\beta_7^{(3,5)}(a_1,a_2)&=-{\frac {143}{5}}\,{a_1}^{2}a_2 ^{4}\left( 12\,{a_1a_2}-5\,a_1-2a_2
 \right), \\
 \beta_6^{(3,5)}(a_1,a_2)&={\frac {11}{5}}\,{a_1}\,{a_2}^{3} \left(  -140\,{a_1}^2\,{a_2}+35\,{a_1}^2+30\,{a_1}\,{a_2}+5\,{a_2}^2-56\,{a_1}\,{a_2}^2+
126\,{a_1}^2\,{a_2}^2\right),\\
\beta_5^{(3,5)}(a_1,a_2)&={a_2}^{2}( {a_2}^3+42\,{a_1}^2{a_2}+28\,{a_1}^3+84\,{a_1}^2\,{a_2}^3-84\,{a_1}^3\,{a_2}^3\\
&+210a_1^3a_2^2-21a_1a_2^3-147a_1^3a_2+15a_1a_2^2-126a_1^2a_2^2 ),\\
\end{align*}
\begin{align*}
\beta_4^{(3,5)}(a_1,a_2)&={\frac {1}{3}}\,{a_2} (
  245 a_1^3  a_2^2  + 21 a_1^3  a_2^4  - 140 a_1^3  a_2 - 140 a_1^3  a_2^3  + 35 a_1 a_2^4
- 75 a_1 a_2^3  - 56 a_1^2  a_2^4  \\
&+ 210 a_1^2  a_2^3  - 210 a_1^2  a_2^2  + 5 a_2^3
         + 56 a_1^2  a_2 + 21 a_1^3  - 5 a_2^4  + 35 a_1 a_2 ) ,\\
\beta_3^{(3,5)}(a_1,a_2)&=a_1^3+{5\over 3}a_1^3a_2^4-{5\over 3}a_1a_2^5-{35\over 3}a_1^3a_2^3+{5\over 3}a_2^3-{50\over 3}a_1a_2^3+{75\over 7}a_1a_2^4
         +{20\over 3}a_1a_2^2\\
&-{50\over 21}a_2^4-10a_1^2a_2^4+6a_1^2a_2+30a_1^2a_2^3
         -{80\over 3}a_1^2a_2^2
         +20a_1^3a_2^2-10a_1^3a_2\\
         &+{5\over 7}a_2^5+{2\over 3}a_1^2a_2^5,\\
    \beta_2^{(3,5)}(a_1,a_2)&= -{1\over 105}(a_1+a_2-1)(140a_1^2a_2^2+126a_1^2-315a_1^2a_2-385a_1a_2^2
     +315a_1a_2\\
     &+70a_1a_2^3-70a_2^3+140a_2^2+5a_2^4),\\
    \beta_1^{(3,5)}(a_1,a_2)&= {1\over 15}(a_1+a_2-1)(3a_1^2+15a_1a_2-15a_1-15a_2+5a_2^2),\\
    \beta_0^{(3,5)}(a_1,a_2)&=-a_1-a_2+1.
\end{align*}
For $n,m\geq 1$ and $a_1=a,\ a_2=1-a$, Berg and Vignat \cite{BV1}  have proved
the following recurrence relation for $\beta^{(n,m)}_{k}(a,1-a)$ which they denoted by $\beta^{(n,m)}_{k}(a)$
\cite[Lemma 3.6]{BV1}:
\begin{equation}\label{relationbeta}
\frac{1}{2k+1}\beta^{(n,m)}_{k+1}(a)=\frac{a^2}{2n-1}
\beta^{(n-1,m)}_{k}(a)+\frac{(1-a)^2}{2m-1}\beta^{(n,m-1)}_{k}(a),
\end{equation}
for $k=0,1,...,m+n-1$. From \eqref{relationbeta} they
derived the positivity of $\beta^{(n,m)}_{k}(a)$  when  $0\leq a\leq 1$  and
also  that $\beta^{(n,m)}_{k}(a)=0$ for $k<\min(m,n)$. Recently, with J. Zeng \cite{az}, we improved this result by giving the explicit single-sum formula for $ \beta_{k}^{(n,m)}(a)$ which was missing in their paper \cite{BV1}.\\
In this paper, our main result is twofold:\\
- for any $a_1,a_2$, a recurrence relation for $\beta^{(n,m)}_{k}(a_1,a_2)$ is given. This recurrence relation reduces to the recurrence system \eqref{relationbeta} when $a_1=a$ and $a_2=1-a$.\\
- for any $a_1,a_2$, an explicit single sum formula for $\beta^{(n,m)}_{k}(a_1,a_2)$, which provides actually the unique solution of the recurrence relation and, then, becomes a generalization of $\beta^{(n,m)}_{k}(a)$ given by Atia and Zeng in \cite{az} when $a_1=a$ and $a_2=1-a$.\\
\begin{lem}
For $n,m\geq 1$, the recurrence relation fulfilled by $\beta_{k}^{(n,m)}(a_1,a_2),\ 0\leq k\leq n+m$ is given by
\begin{equation}\label{eq:equ1}
\beta_{n+m}^{(n+1,m-1)}(a_1,a_2)-{a_1^2(2m-1)(2m+1)\over a_2^2(2n-1)(2n+1)}\beta_{n+m}^{(n-1,m+1)}(a_1,a_2)=0,
\end{equation}
and for $0\leq k\leq n+m-1$, we have
\begin{eqnarray}\label{eq:equ2}
  \beta_{k}^{(n+1,m-1)}(a_1,a_2)-{a_1^2(2m-1)(2m+1)\over a_2^2(2n-1)(2n+1)}\beta_{k}^{(n-1,m+1)}(a_1,a_2)&&  \nonumber\\
  =\beta_{k}^{(n,m-1)}(a_1,a_2)-{a_1^2(2m-1)
(2m+1)\over a_2^2(2n-1)(2n+1)}\beta_{k}^{(n-1,m)}(a_1,a_2).  & &
\end{eqnarray}
\end{lem}
{\bf Proof.} In one hand we have
$$q_{n+1}(a_1u)q_{m-1}(a_2u)=\sum_{k=0}^{n+m}\beta_{k}^{(n+1,m-1)}(a_1,a_2)q_{k}(u),$$
in the other hand, using \eqref{recu_bess}, we have,
\begin{eqnarray*}
 && q_{n+1}(a_1u)q_{m-1}(a_2u) = \biggl (q_{n}(a_1u)+ {a_1^2u^2\over (2n-1)(2n+1)}q_{n-1}(a_1u) \biggr )q_{m-1}(a_2u) \\
   &=& q_{n}(a_1u)q_{m-1}(a_2u) + {a_1^2u^2\over (2n-1)(2n+1)}q_{n-1}(a_1u)q_{m-1}(a_2u) \\
   &=& q_{n}(a_1u)q_{m-1}(a_2u)  + {a_1^2(2m-1)(2m+1)\over a_2^2(2n-1)(2n+1)}q_{n-1}(a_1u){a_2^2u^2\over (2m-1)(2m+1)}q_{m-1}(a_2u) \\
&=& q_{n}(a_1u)q_{m-1}(a_2u) + {a_1^2(2m-1)(2m+1)\over a_2^2(2n-1)(2n+1)}q_{n-1}(a_1u)\biggl (q_{m+1}(a_2u)-q_{m}(a_2u)\biggr )\nonumber
\end{eqnarray*}
where we used again \eqref{recu_bess}, finally, we obtain
$$
  q_{n+1}(a_1u)q_{m-1}(a_2u)-{a_1^2(2m-1)(2m+1)\over a_2^2(2n-1)(2n+1)}q_{n-1}(a_1u)q_{m+1}(a_2u) $$
  $$= q_{n}(a_1u)q_{m-1}(a_2u) - {a_1^2(2m-1)(2m+1)\over a_2^2(2n-1)(2n+1)}q_{n-1}(a_1u)q_{m}(a_2u),$$
and because of the degree of polynomials $q_k(u)$ we have \eqref{eq:equ1}
and for $0\leq k\leq n+m-1$ we have \eqref{eq:equ2}.
\bigskip

\begin{thm}\label{thm:A-Z}
For $i=0,1,...,n+m$, we have
\begin{align}\label{eq:thm}
&\beta_{k}^{(n,m)}(a_1,a_2)={a_1^{-m+k}a_2^{-n+k}
(1/2)_k\over 4^{m+n-k}(m+n-k)!(1/2)_n(1/2)_m}\nonumber \\
  &\sum_{i=0}^{m+n-k}{a_1^{m+n-k-i}(^{m+n-k}_i)
(n+1-i)_{2i}}\nonumber\\
&\sum_{j=0}^{m+n-k-i}{(-1)^j(^{m+n-k-i}_j)
(-n+k+j+i+1)_{2(m+n-k-i-j)}(k+2-j)_{2j}
a_2^{j+i}}.
\end{align}
\end{thm}
which we write using $_3F_2$ hypergeometric functions as
\begin{thm}\label{thm:A-J}
For $i=0,1,...,n+m$, we have
\begin{eqnarray}\label{eq:hyp}
  \beta_{k}^{(n,m)}(a_1,a_2)&=& {a_1^{-m+k}a_2^{m}(1/2)_k
  \over 4^{m+n-k}(m+n-k)!(1/2)_n(1/2)_m}\nonumber   \\
  &&\sum_{i=0}^{m+n-k}{a_1^{i}(^{m+n-k}_{m+n-k-i})
(-m+k+i+1)_{2(m+n-k-i)}(m-i+1)_{2i}} \nonumber \\
 && \ _3F_2(\begin{array}{cc}
        k+2,\ -k-1,\ -i \\
        -m-i,\ m-i+1
      \end{array};a_2)a_2^{-i}.
\end{eqnarray}
\end{thm}
\bigskip
\noindent
{\bf Remarks.}\\
1. This formula was deduced using the same approach done in \cite{az} pages 4 and 5 by, just, changing $a$ by $a_1$ and $1-a$ by $a_2$.\\
 2. To compute this formula with, for example, Maple, one should compute $\beta_{n+m}^{(n,m)}(a_1,a_2)$, $\beta_{n+m-1}^{(n,m)}(a_1,a_2)$,..., $\beta_{0}^{(n,m)}(a_1,a_2)$ and then replace $n,m$ by their values (please see the Maple program given in the end of this paper).\\

\bigskip
\noindent
{\bf Proof of theorem 3.} Let us, first, prove that \eqref{eq:hyp} fulfils \eqref{eq:equ1}:
$$\beta_{n+m}^{(n+1,m-1)}(a_1,a_2)={a_1^{n+1}a_2^{m-1}\sqrt\pi\Gamma(1/2+n+m)\over \Gamma(n+3/2)\Gamma(m-1/2)},$$
and
$$\beta_{n+m}^{(n-1,m+1)}(a_1,a_2)={a_1^{n-1}a_2^{m+1}\sqrt\pi\Gamma(1/2+n+m)\over \Gamma(n-1/2)\Gamma(m+3/2)},$$
then
$${a_1^2(2m-1)(2m+1)\over a_2^2(2n-1)(2n+1)}\beta_{n+m}^{(n-1,m+1)}(a_1,a_2)={a_1^{n+1}a_2^{m-1}\sqrt\pi\Gamma(1/2+n+m)\over \Gamma(n+3/2)\Gamma(m-1/2)}.$$

Second, we prove that \eqref{eq:hyp} fulfils \eqref{eq:equ2}, so let us substract the rhs from lhs of \eqref{eq:equ2} to obtain
$$\beta_{k}^{(n+1,m-1)}(a_1,a_2)-{a_1^2(2m-1)(2m+1)\over a_2^2(2n-1)(2n+1)}\beta_{k}^{(n-1,m+1)}(a_1,a_2)$$

$$-\beta_{k}^{(n,m-1)}(a_1,a_2)+{a_1^2(2m-1)(2m+1)\over a_2^2(2n-1)(2n+1)}\beta_{k}^{(n-1,m)}(a_1,a_2)$$

$$=
{a_1^{-m+1+k}a_2^{m-1}\sqrt\pi\Gamma(1/2+k)\over 4^{m+n-k}(m+n-k)!\Gamma(n+3/2)\Gamma(m-1/2)}\times$$
$$\biggl (\sum_{i=0}^{m+n-k}{a_1^{i}(^{m+n-k}_{m+n-k-i})(-m+k+i+2)_{2(n+m-k-i)}
(m-i)_{2i}}$$
$$_3F_2([-i, k+2, -k-1],[-m+1-i, m-i],a_2)a_2^{-i}
$$
$$-\sum_{i=0}^{m+n-k}{a_1^{i}(^{m+n-k}_{m+n-k-i})
(-m+k+i)_{2(n+m-k-i)}(m-i+2)_{2i}}\nonumber$$
$$_3F_2([-i, k+2, -k-1],[-m-i-1, m-i+2],a_2)a_2^{-i}\biggr )
$$

$$
-{a_1^{-m+1+k}a_2^{m-1}\sqrt\pi\Gamma(1/2+k)\over 4^{m+n-k-1}(m+n-k-1)!\Gamma(n+1/2)\Gamma(m-1/2)}\times\nonumber $$
$$\biggl (\sum_{i=0}^{m+n-k-1}a_1^{i}(^{m+n-k-1}_{m+n-k-i-1})
(-m+k+i+2)_{2(n+m-k-i-1)}(m-i)_{2i}\nonumber$$
$$_3F_2([-i, k+2, -k-1],[-m-i+1, m-i],a_2)a_2^{-i}
$$
$$
-{(m+1/2)a_1\over (n+1/2)a_2}\sum_{i=0}^{m+n-k-1}{a_1^{i}(^{m+n-k-1}_{m+n-k-i-1})
(-m+k+i+1)_{2(n+m-k-i-1)}(m-i+1)_{2i}}
$$
$$_3F_2([-i, k+2, -k-1],[-m-i, m-i+1],a_2)a_2^{-i}\biggr ),$$
because

$${(2m-1)(2m+1)\over (2n-1)(2n+1)\Gamma(n-1/2)\Gamma(m+3/2)}
={1\over \Gamma(n+3/2)\Gamma(m-1/2)}.$$
Cancelling the common factor
$${a_1^{-m+1+k}a_2^{m-1}\sqrt\pi\Gamma(1/2+k)\over 4^{m+n-k-1}(m+n-k-1)!\Gamma(n+1/2)\Gamma(m-1/2)}$$ in both quantities, we get
$$
{1\over 4(m+n-k)(n+1/2)}\times\nonumber$$
$$\biggl (\sum_{i=0}^{m+n-k}{a_1^{i}(^{m+n-k}_{m+n-k-i})(-m+k+i+2)_{2(n+m-k-i)}
(m-i)_{2i}}\nonumber$$
$$_3F_2([-i, k+2, -k-1],[-m+1-i, m-i],a_2)a_2^{-i}
$$
$$-\sum_{i=0}^{m+n-k}{a_1^{i}(^{m+n-k}_{m+n-k-i})
(-m+k+i)_{2(n+m-k-i)}(m-i+2)_{2i}}\nonumber$$
$$_3F_2([-i, k+2, -k-1],[-m-i-1, m-i+2],a_2)a_2^{-i}\biggr )
$$

$$-\biggl (\sum_{i=0}^{m+n-k-1}a_1^{i}(^{m+n-k-1}_{m+n-k-i-1})
(-m+k+i+2)_{2(n+m-k-i-1)}(m-i)_{2i}\nonumber$$
$$_3F_2([-i, k+2, -k-1],[-m-i+1, m-i],a_2)a_2^{-i}
$$
$$
-{(m+1/2)a_1\over (n+1/2)a_2}\sum_{i=0}^{m+n-k-1}{a_1^{i}(^{m+n-k-1}_{m+n-k-i-1})
(-m+k+i+1)_{2(n+m-k-i-1)}(m-i+1)_{2i}}
$$
$$_3F_2([-i, k+2, -k-1],[-m-i, m-i+1],a_2)a_2^{-i}\biggr ),$$
equivalently
$$
{1\over 4(m+n-k)(n+1/2)}\times\nonumber$$
$$\biggl(\sum_{i=0}^{m+n-k}{a_1^{i}(^{m+n-k}_{m+n-k-i})(-m+k+i+2)_{2(n+m-k-i)}
(m-i)_{2i}}\nonumber$$
$$_3F_2([-i, k+2, -k-1],[-m+1-i, m-i],a_2)a_2^{-i}
$$
$$-\sum_{i=0}^{m+n-k}{a_1^{i}(^{m+n-k}_{m+n-k-i})
(-m+k+i)_{2(n+m-k-i)}(m-i+2)_{2i}}\nonumber$$
$$_3F_2([-i, k+2, -k-1],[-m-i-1, m-i+2],a_2)a_2^{-i}\biggr )
$$

$$-\sum_{i=0}^{m+n-k-1}a_1^{i}(^{m+n-k-1}_{m+n-k-i-1})
(-m+k+i+2)_{2(n+m-k-i-1)}(m-i)_{2i}\nonumber$$
$$_3F_2([-i, k+2, -k-1],[-m-i+1, m-i],a_2)a_2^{-i}
$$
$$
+{(m+1/2)a_1\over (n+1/2)a_2}\sum_{i=0}^{m+n-k-1}{a_1^{i}(^{m+n-k-1}_{m+n-k-i-1})
(-m+k+i+1)_{2(n+m-k-i-1)}(m-i+1)_{2i}}
$$
$$_3F_2([-i, k+2, -k-1],[-m-i, m-i+1],a_2)a_2^{-i}.$$
To prove that this expression vanishes, it suffices to prove that the coefficient of $a_1^j$ vanishes. The coefficient of $a_1^j$ is given by

$$
{1\over 4(m+n-k)(n+1/2)}\times\nonumber$$
$$\biggl({a_1^{j}(^{m+n-k}_{m+n-k-j})(-m+k+j+2)_{2(n+m-k-j)}
(m-j)_{2j}}\nonumber$$
$$_3F_2([-j, k+2, -k-1],[-m+1-j, m-j],a_2)a_2^{-j}
$$
$$-{a_1^{j}(^{m+n-k}_{m+n-k-j})
(-m+k+j)_{2(n+m-k-j)}(m-j+2)_{2j}}\nonumber$$
$$_3F_2([-j, k+2, -k-1],[-m-j-1, m-j+2],a_2)a_2^{-j}\biggr )
$$

$$-a_1^{j}(^{m+n-k-1}_{m+n-k-j-1})
(-m+k+j+2)_{2(n+m-k-j-1)}(m-j)_{2j}\nonumber$$
$$_3F_2([-j, k+2, -k-1],[-m-j+1, m-j],a_2)a_2^{-j}
$$
$$
+{(m+1/2)a_1\over (n+1/2)a_2}{a_1^{j-1}(^{m+n-k-1}_{m+n-k-j})
(-m+k+j)_{2(n+m-k-j)}(m-j+2)_{2(j-1)}}
$$
$$_3F_2([-(j-1), k+2, -k-1],[-m-(j-1), m-(j-1)+1],a_2)a_2^{-(j-1)}.$$\medskip

A short computation (with Maple) of this quantity gives zero: \\
$Q1:=((a1^i*binomial(n+m-k,n+m-k-i)*pochhammer(-m+k+i+2,2*n+2*m-2*k-2*i)*pochhammer(m-i,2*i)
*hypergeom([-i, k+2, -k-1],[-m-i+1, m-i],a2)*a2^{(-i)})-(a1^i*binomial(n+m-k,n+m-k-i)*pochhammer(-m+k+i,2*n+2*m-2*k-2*i)
*pochhammer(m-i+2,2*i)*hypergeom([-i, k+2, -k-1],[m-i+2, -m-i-1],a2)*a2^{(-i)}))
-4*(n+m-k)*(n+1/2)*(a1^i*binomial(n+m-k-1,n+m-k-i-1)*pochhammer(-m+k+i+2,2*n+2*m-2*k-2*i-2)
*pochhammer(m-i,2*i)*hypergeom([-i, k+2, -k-1],[-m-i+1, m-i],a2)*a2^{(-i)});
Q2:=-4*(n+m-k)*(m+1/2)*a1/a2*(a1^{(i-1)}*binomial(n+m-k-1,n+m-k-(i-1)-1)
*pochhammer(-m+k+(i-1)+1,2*n+2*m-2*k-2*(i-1)-2)*pochhammer(m-(i-1)+1,2*(i-1))
*hypergeom([-(i-1), k+2, -k-1],[-m-(i-1), m-(i-1)+1],a2)*a2^{(-i+1)});
simplify(Q1-Q2);$ \\

{\bf Particular case.}\\
Let us prove that \eqref{eq:equ2} reduces to \eqref{relationbeta} when $a_1=a,\ a_2=1-a$.\\
From \eqref{eq:equ2} we have
\begin{eqnarray}\label{relationbeta2}
{(1-a)^2\over 2m-1}\beta_{k}^{(n,m-1)}(a,1-a)-{a^2
(2m+1)\over (4n^2-1)}\beta_{k}^{(n-1,m)}(a,1-a)
&&  \nonumber\\
  =  {(1-a)^2\over 2m-1}\beta_{k}^{(n+1,m-1)}(a,1-a)-{a^2(2m+1)\over (4n^2-1)}\beta_{k}^{(n-1,m+1)}(a,1-a).  & &
\end{eqnarray}
equivalently
\begin{eqnarray}\label{relationbeta3}
{(1-a)^2\over 2m-1}\beta_{k}^{(n,m-1)}(a,1-a)={a^2
(2m+1)\over (4n^2-1)}\beta_{k}^{(n-1,m)}(a,1-a)
&&  \nonumber\\
  +  {(1-a)^2\over 2m-1}\beta_{k}^{(n+1,m-1)}(a,1-a)-{a^2(2m+1)\over (4n^2-1)}\beta_{k}^{(n-1,m+1)}(a,1-a).  & &
\end{eqnarray}
Adding $\frac{a^2}{2n-1}
\beta^{(n-1,m)}_{k}(a,1-a)$ to both sides, we get
\begin{eqnarray}\label{relationbeta3}
{(1-a)^2\over 2m-1}\beta_{k}^{(n,m-1)}(a,1-a)+\frac{a^2}{2n-1}
\beta^{(n-1,m)}_{k}(a,1-a)&&  \nonumber\\=
{a^2
(2m+1)\over (4n^2-1)}\beta_{k}^{(n-1,m)}(a,1-a)
  +  {(1-a)^2\over 2m-1}\beta_{k}^{(n+1,m-1)}(a,1-a)&&  \nonumber\\
  -{a^2(2m+1)\over (4n^2-1)}\beta_{k}^{(n-1,m+1)}(a,1-a)+\frac{a^2}{2n-1}
\beta^{(n-1,m)}_{k}(a,1-a).  & &
\end{eqnarray}
According to \eqref{relationbeta} le lhs is equal to $\frac{1}{2k+1}\beta^{(n,m)}_{k+1}(a,1-a)$.\\
Using \eqref{eq:lin}, the rhs becomes
\begin{eqnarray}\label{relationbeta5}
 - {a^2(2m+1)\over (4n^2-1)}q_{n-1}(au)q_{m+1}((1-a)u)+\frac{a^2}{2n-1}(1+{2m+1\over 2n+1})q_{n-1}(au)q_{m}((1-a)u)&& \nonumber\\
  +{(1-a)^2\over 2m-1}q_{n+1}(au)q_{m-1}((1-a)u).&&\nonumber
\end{eqnarray}
Using \eqref{recu_bess}, we obtain
\begin{eqnarray}\label{relationbeta5}
 - {a^2(2m+1)\over (4n^2-1)}q_{n-1}(au)\biggl(q_m((1-a)u)+{(1-a)^2u^2\over 4m^2-1}q_{m-1}((1-a)u)\biggr)&&\nonumber\\
 +\frac{a^2}{2n-1}(1+{2m+1\over 2n+1})q_{n-1}(au)q_{m}((1-a)u)&&\nonumber \\
  +{(1-a)^2\over 2m-1}\biggl (q_n(au)+{a^2u^2\over 4n^2-1}q_{n-1}(au)\biggr )q_{m-1}((1-a)u).&&\nonumber
\end{eqnarray}
After simplification, we get the rhs of \eqref{relationbeta}.

\section {Applications}

1- These coefficients $\beta_{k}^{(n,m)}(a_1,a_2)$ with $a_1+a_2 \ne 1$ can be applied in: if $X$ and $Y $ are two
 student random variables with $n$ and $m$ degrees of freedom then
the linear combination $a_1X + a_2Y$ has for characteristic function
$$e^{(-a_1u-a_2u)}q_n(a_1u)q_m(a_2u)= e^{(-a_1u-a_2u)} \sum_{k=0}^{n+m} \beta_k^{n,m}(a_1,a_2) q_k(u),$$
On the other hand, we have
$$a_1X+a_2Y = {(a_1+a_2)}({a_1\over {(a_1+a_2)}}X + {a_2\over {(a_1+a_2)}}Y) = {(a_1+a_2)} (\tilde{a_1} X+ \tilde{a_2} Y)$$
with $\tilde{a_1} + \tilde{a_2}=1$
then it exists a NON TRIVIAL relation between the coefficients $\beta_{k}^{(n,m)}(a_1,a_2)$ and the coefficients $\beta(a_1,1-a_1)$ which is not clear in their expressions.\\
2- For $a_1=a,\ a_2=1-a$, these coefficients $\beta_{k}^{(n,m)}(a,1-a)$
give a formula reducing a sum of hypergeometric functions $_3F_2$ to $_2F_1$:
\begin{thm}\label{eq:thm32} Taking into account \eqref{eq:hyp} and formulas $(7)-(8)$ given in \cite{az}, we get: for $k\geq \lceil( n+m-1)/2\rceil$
\begin{equation*}
{a^{2n+2m-2k}(1-a)^{-m-n+k}\Gamma(n+m+2)\over\Gamma(-n-m+2k+2)}
\ _2F_1(\begin{array}{c}
        -m+k+1,\ -2n-2m+2k \\
        -n-m+2k+2
      \end{array};{1\over a})
\end{equation*}
\begin{eqnarray*}
=\sum_{i=0}^{m+n-k}{a^{i}(^{m+n-k}_{m+n-k-i})
(-m+k+i+1)_{2(m+n-k-i)}(m-i+1)_{2i}}\\
_3F_2(\begin{array}{c}
        k+2,\ -k-1,\ -i \\
        -m-i,\ m-i+1
      \end{array};1-a
)(1-a)^{-i}
\end{eqnarray*}
and for $k\leq \lfloor (n+m-1)/2\rfloor $
\begin{equation*}
{(-a)^{n+1+m}(1-a)^{-m-n+k}\Gamma(2n+2m-2k+1)\Gamma(n-k)\over\Gamma(n+m-2k)\Gamma(-m+k+1)}
\ _2F_1(\begin{array}{c}
        n-k,\ -n-m-1 \\
        n+m-2k
      \end{array};{1\over a})
\end{equation*}
\begin{eqnarray*}
=\sum_{i=0}^{m+n-k}{a^{i}(^{m+n-k}_{m+n-k-i})
(-m+k+i+1)_{2(m+n-k-i)}(m-i+1)_{2i}}\\
_3F_2(\begin{array}{c}
        k+2,\ -k-1,\ -i \\
        -m-i,\ m-i+1
      \end{array};1-a
)(1-a)^{-i}
\end{eqnarray*}
\end{thm}

{\noindent\bf Acknowledgement. } This work was supported by the research unit ur11es87, Gabes university, Tunisia. I would like to thank C. Vignat for pointing out the first application of these coefficients.

\bigskip
Please find next a Maple program which, not only, tests that our formula is right from $min(n,m)$ to $n+m$ but, also, show that\\  \noindent
$\beta_k^{n,m}(a,1-a)=0$  for $k<min(n,m)$.

\noindent
${>restart;}\\
>A:= (n,m)->q(n,a1*u)*q(m,a2*u)-sum( beta(n,m,k,a1,a2) \\
   \qquad  *q(k,u),k=min(n,m)..n+m :)$\\
We assume n less or equal m. This program runs from $min(n,m)$ untill $n+m$, take
any values of $n,\ m$, for example 2 and 8
\\
${>AA:=A(2,8):;}$\\

${>alpha:=(n,k)->n!*(2*n-k)!*2^k/(2*n)!/(n-k)!/k!:}$\\
${>q:=(n,u)->sum(alpha(n,k)*u^k,k=0..n):;}$\\
\begin{eqnarray*}
  >beta:=&&(n,m,k,a1,a2)->factor(a1^{(-m+k)}*a2^m*Pi^{(1/2)}  \\
   &&*GAMMA(1/2+k)*
sum(a1^i*binomial(n+m-k,n+m-k-i)\\
&&*pochhammer(-m+k+i+1,2*n+2*m-2*k-2*i)\\
&&*pochhammer(m-i+1,2*i)*simplify(hypergeomt(n,m,i,k))\\
&&*a2^{(-i)},i = 0 ..
n+m-k)/(4^(n+m-k))/(n+m-k)!\\
&&/GAMMA(n+1/2)/GAMMA(m+1/2)):; 
\end{eqnarray*}
${>AAA:=factor(AA):;}$\\
\begin{eqnarray*}
   &&>hypergeomt:=(n,m,i,k)->simplify(hypergeom([-i, k+2, -k-1],  \\
   && [-m-i,m-i+1],a2)): 
\end{eqnarray*}
${>collect(factor(simplify(AAA)),u);}$\\
$$ 1/5*(-1+a1+a2)*(5*a2*a1-5*a1-5*a2+2*a2^2)*u $$
$$ +1/5*(-1+a1+a2)*(5*a2*a1-5*a1-5*a2+2*a2^2-5);$$
We meet again that $beta(n,m,k,a,1-a)$ vanish for
$k<min(n,m)$.\eject

\end{document}